\documentclass[a4paper,14pt]{article}
\usepackage[14pt]{extsizes}
\usepackage{mathtext}                    % Русский язык в формулах
\usepackage{cmap}                        % Поддержка поиска русских слов в PDF (pdflatex)
\usepackage[cp1251]{inputenc}            % Выбор языка и кодировки
\usepackage[english, russian]{babel}     % Настройки для русского языка
\usepackage[left=2.5cm,right=1cm,top=1cm,bottom=2cm]{geometry} % поля страницы
\usepackage{amssymb}
\usepackage{amsmath}

\newtheorem{definition}{Определение}
\newtheorem{corollary}{Следствие}
\newtheorem{remark}{Замечание}
\newtheorem{lemma}{Лемма}
\newtheorem{theorem}{Теорема}

\begin{document}

\centerline{\textbf{Некоторые вопросы теории возмущений в функциональном исчислении}}
\centerline{\textbf{замкнутых операторов в банаховом пространстве}}

\centerline{\textbf{А. Р. Миротин}}

\centerline{{amirotin@yandex.ru}}

The work is a continuation of research of A.A. Atvinovskii and the author of functional calculus of closed operators on Banach spaces based on Markov and related functions as symbols. The following topics in the perturbation theory are considered:  estimates of bounded perturbations of  operator functions with respect to general operator ideal norms, operator Lipschitzness, moment inequality, Freshet  operator differentiability, analyticity of operator functions under consideration with respect to the  perturbation parameter, spectral shift function and Livshits–Krein trace formulae.

\section{Введение и предварительные сведения} \label{ }

Различным проблемам теории возмущений линейных операторов, относящимся к классическому функциональному исчислению фон Неймана-Мурье-Данфорда самосопряженных, унитарных и нормальных операторов в гильбертовом пространстве, посвящена обширная литература (см., например, \cite{Peller} --- \cite{Kissin_2012}). Данная работа посвящена аналогичным вопросам теории возмущений линейных операторов, возникающим в функциональном исчислении замкнутых операторов в банаховом пространстве, построенном в \cite{mirotin_AM},  \cite{mirotin_AM_2} (см. также   \cite{AMpfmt2013}, \cite{Obr_resolvent}). Вместе с тем, ее можно читать независимо от указанных работ.

В статье  рассматриваются следующие вопросы теории возмущений:  оценки ограниченных возмущений значений   рассматриваемых в работе операторных функций в нормах общих операторных идеалов, операторная и коммутаторная липшицевость таких функций, неравенство моментов, операторная дифференцируемость по Фреше, аналитическая зависимость значений   рассматриваемых операторных функций от параметра возмущения, функция спектрального сдвига и аналог формулы следа Лифшица-Крейна для ядерных возмущений.

Единственной известной автору работой по операторно и коммутаторно  липшицевым функциям в контексте банаховых пространств является недавняя статья \cite{RST}\footnote{Автор благодарит рецензента, указавшего ему на эту работу.}. В этой работе рассмотрено функциональное исчисление для операторов скалярного типа в банаховых пространствах, использующее в качестве символов довольно широкий класс функций, тесно связанный с пространством Соболева $W^{1,2}(\mathbb{R})$ и однородным пространством Бесова $\dot B^{1}_{\infty,1}(\mathbb{R})$. С помощью развитой  теории двойных операторных интегралов для операторов скалярного типа в банаховых пространствах в указанной статье получены, в частности,  оценки  коммутаторов весьма общего вида в операторной норме и нормах операторных идеалов, удовлетворяющих некоторым естественным ограничениям.

Напомним необходимые понятия и результаты из \cite{mirotin_AM} и \cite{mirotin_AM_2}.

\begin{definition}
 Пусть $a<b$. Говорят, что функция $g$ относится к классу Маркова $R[a,b]$, если она голоморфна в верхней полуплоскости, отображает ее в свое замыкание, а также голоморфна и положительна на $(-\infty,a)$ и голоморфна и отрицательна на $(b,\infty)$.
\end{definition}

Скажем, что функция $g$ принадлежит классу $R^0$, если она принадлежит $R[0,b]$ и непрерывна в нуле.

 \begin{definition} Обозначим через $ZR(0,b]\ (ZR[a,b])$ множество функций вида $zg(z)$,  где $g\in R^0$ (соответственно, $g\in R[a,b]$).
 \end{definition}

В силу  интегрального представления функций Маркова  \cite[теорема П.6]{krein}, функции из $ZR(0,b]$ имеют в точности  вид
$$
f(z)=\int\limits_0^b\frac{z}{t-z}d\tau(t),\eqno(1)
$$
где $\tau$ --- единственная ограниченная положительная мера на  $(0,b]$, такая, что $\int_0^bd\tau(t)/t<\infty$ ("представляющая мера"). Аналогичное интегральное представление имеют и функции из $ZR[a,b]$.

Следующая лемма дает внутреннюю характеристику функций из $ZR(0,b]$.

 \begin{lemma}
 Функция $f$ принадлежит  $ZR(0,b]$ тогда и только тогда, когда она голоморфна и отрицательна на промежутках $(-\infty ,0)$ и $(b,+\infty )$, причем функция  $f(z)/z$ непрерывна в нуле и отображает открытую верхнюю полуплоскость $\{z\in \mathbb{C}: {\rm Im}z>0\}$ в свое замыкание.
\end{lemma}

Доказательство.
Достаточно заметить, что указанные свойства функции $f$ равносильны тому, что функция  $g(z):=f(z)/z$ принадлежит классу $R[0,b]$ функций Маркова и непрерывна в нуле, т. е. принадлежит классу $R^0$ в смысле \cite{mirotin_AM_2}.
%\end{proof}

Аналогичным образом можно дать внутреннюю характеристику функций из $ZR[0,b]$.

\textbf{Примеры 1}. Функция
$$
\frac{z}{b-z}\left(\frac{z}{z-b}\right)^{\alpha-1}=\frac{\sin\pi\alpha}{\pi}\int\limits_0^b\frac{z}{t-z}t^{\alpha-1}(b-t)^{-\alpha}dt \quad (0<\alpha<1)
$$
принадлежит  $ZR[0,b]$, но не принадлежит  $ZR(0,b]$;

функция
$$
 z\left(1-\left(\frac{z}{z-b}\right)^\alpha\right)=\frac{\sin\pi\alpha}{\pi}\int\limits_0^b\frac{z}{t-z}\left(\frac{t}{b-t}\right)^\alpha dt \quad (0<\alpha<1)
$$
принадлежит  $ZR(0,b]$ (см. \cite[2.2.6, формулы № 7, 5]{PBM}).

Всюду ниже $A$ --- замкнутый плотно определенный оператор в комплексном банаховом пространстве $X$,  $\sigma(A)$ и  $\rho(A)$ --- спектр и резольвентное множество оператора $A$ соответственно, $LB(X)$ --- алгебра линейных ограниченных операторов в  $X$.

\begin{definition}
Будем говорить, что (замкнутый плотно определенный) оператор $A$  в пространстве $X$ принадлежит классу $V_{(0,b]}(X)$, если  $(0,b]\subset \rho(A)$ и для некоторой постоянной $M_A>0$  выполняется неравенство
$$
\|R(t,A)\| \leq \frac{M_A}{t}, t\in (0,b].
$$
\end{definition}

Класс $V_{(0,b]}(X)$ весьма широк. Например, ему принадлежат  все  операторы вида  $-A$, где $A$  есть слабо позитивный оператор \cite{Pys}. Укажем еще один класс операторов из  $V_{(0,b]}(X)$. Напомним, что ограниченный оператор $T$ в пространстве $X$ называется \textit{оператором Ритта}, если $\sigma(T)\subset \overline{\mathbb{D}}$ и
$\|R(t,T)\|\leq C/|t-1|\ (|t|>1)$ (см., например, \cite{1}, \cite{2}).

\textbf{Пример 2}. Если  оператор Ритта $T$ в пространстве $X$ имеет левый обратный  $T^{-1}$, то    $T^{-1}+I\in V_{(0,1]}(X)$. В самом деле, если $1/t\in\rho(T)$ (в частности, если $|t|<1$), то легко проверить, что оператор
$(-1/t)TR(1/t,T)=(-1/t)R(1/t,T)T$   равен $R(t,T^{-1})$. Следовательно, $(0,1]\subset \rho(T^{-1}+I)$
и при $t\in (0,1]$ имеем
$$
\|R(t,T^{-1}+I)\|=\|R(t-1,T^{-1})\|\leq\|1/(t-1)TR(1/(t-1),T)\|\leq C\|T\|/t.
$$

 \begin{definition} Для любой функции $f\in ZR(0,b]$ с представляющей мерой $\tau$ и любого $A\in V_{(0,b]}(X)$ положим
 $$
 f(A)=\int\limits_0^bAR(t,A)d\tau(t). \eqno(2)
$$
\end{definition}

Из хорошо известного тождества
$$
AR(t,A)=-I+tR(t,A)\quad (t\in\rho(A))\eqno(3)
$$
следует, что при $A\in V_{(0,b]}(X)$ интеграл в (2) существует в смысле Бохнера (относительно операторной нормы), оператор  $f(A)$ определен на $X$ и ограничен, причем $\|f(A)\|\leq (1+M_A)\tau([0,b])$.

Ясно, что  $f(A)=Ag(A)$, где $g(A)$  понимается в смысле  \cite{mirotin_AM_2}, т. е. $g(A)=\int_0^bR(t,A)d\tau(t)$. Поэтому данное исследование можно рассматривать как продолжение работ  \cite{mirotin_AM} --- \cite{mirotin_AM_2}.

\section{ Оценка возмущений ограниченными операторами}

\begin{theorem}
Пусть $f\in  ZR(0,b]$. Для любых операторов $A, B\in V_{(0,b]}(X)$, таких, что $D(A)\supseteq D(B)$ и оператор $A-B$ ограничен, справедливо неравенство
$$
\|f(A)-f(B)\|\leq -(M_A+M_B+M_AM_B)f(-\|A-B\|).
$$
\end{theorem}

Доказательство. Ясно, что
$$
 f(A)-f(B)=\int\limits_0^b(AR(t,A)-BR(t,B))d\tau(t).\eqno(4)
$$
 Рассмотрим подынтегральную функцию $\varphi(t):=AR(t,A)-BR(t,B)$. Поскольку $D(A)\supseteq D(B)$, то в силу формулы (3) и второго резольвентного тождества имеем  $\varphi(t)=t(R(t,A)-R(t,B))=tR(t,A)(A-B)R(t,B)$, откуда  $\|\varphi(t)\|\leq M_AM_B\|A-B\|/t$. Но так как $\|\varphi(t)\|\leq M_A+M_B$, то
$$
\|\varphi(t)\|\leq M_AM_B\|A-B\|\min\left\{\frac{ M_A+M_B}{M_AM_B\|A-B\|},\frac{1}{t}\right\}.
$$
Полагая в неравенстве
$$
\min\left\{a,\frac{1}{t}\right\}\leq\frac{1+ad}{t+d},
$$
справедливом для всех $a,d,t>0$,  $a=(M_A+M_B)/(M_AM_B\|A-B\|), d=\|A-B\|$, получим
$$
\|\varphi(t)\|\leq \frac{(M_A+M_B+M_AM_B)\|A-B\|}{t+\|A-B\|}.
$$
Теперь  из (4) следует, что
$$
\|f(A)-f(B)\|\leq(M_A+M_B+M_AM_B)\|A-B\|\int\limits_0^b\frac{d\tau(t)}{t+\|A-B\|}=
$$
$$
-(M_A+M_B+M_AM_B)f(-\|A-B\|),
 $$
что и требовалось доказать.
%\end{proof}

 \begin{theorem}
 Если в условиях теоремы 1 операторы $A-B$ и $R(t,B)$ коммутируют, то при всех $x\in D(A)$ справедливо неравенство
$$
\|(f(A)-f(B))x\|\leq -(M_A+M_B+M_AM_B)f(-\|(A-B)x\|).
$$
\end{theorem}

Эта теорема доказывается так же, как и теорема  1.

Следующее неравенство по своей формулировке аналогично   неравенствам моментов, установленным в \cite{Pys}, \cite{SMZ}, и есть частный случай теоремы 2 (при $B=O$).

\begin{corollary}\cite{MA}
 Пусть $A\in V_{(0,b]}(X)$, $f\in ZR(0,b]$. Тогда для любого $x\in D(A), \|x\|=1$ справедливо неравенство
$$
\|f(A)x\|\leq -(2M_A+1)f(-\|Ax\|).
$$
\end{corollary}

\begin{corollary}
Если $A\in V_{(0,b]}(X)$, $f\in ZR(0,b]$, то при всех  $x\in D(A), \|x\|=1$
$$
\|f(A)x\|\leq (2M_A+1)f'(-0)\|Ax\|.
$$
\end{corollary}

Доказательство. Заметим, что функция $f(-s)/(-s)$  убывает на множестве  $s>0$. Поэтому
$$
f'(-0)=\lim_{s\to+0}\frac{f(-s)}{-s}\geq \frac{f(-s)}{-s} \ (s>0).
$$
В частности,  $-f(-\|Ax\|)\leq f'(-0)\|Ax\|$.
%\end{proof}

 Далее $(\mathcal{I},\|\cdot\|_{\mathcal{I}})$ обозначает \textit{операторный идеал} в   $X$, т. е. двусторонний идеал алгебры $LB(X)$ ограниченных операторов в пространстве $X$, полный относительно нормы $\|\cdot\|_{\mathcal{I}}$,  удовлетворяющей условиям  $\|ASB\|_{\mathcal{I}}\leq \|A\|\|S\|_{\mathcal{I}}\|B\|$, $\|S\|\leq \|S\|_{\mathcal{I}}$
 при всех $A, B\in LB(X)$ и $S\in \mathcal{I}$ (случай $\mathcal{I}=LB(X)$  не исключается  и представляет интерес).

 Следующая теорема показывает, что функции класса  $ZR(0,b]$ сохраняют возмущения операторами из    $\mathcal{I}$ и, в частности,  являются  $\mathcal{I}$-\textit{операторно липшицевыми} (см., например, \cite{Peller})  в классе операторов $V_{(0,b]}^c(X):=$ $\{A\in V_{(0,b]}(X): M_A\leq c\}$ $(c={\rm const})$.

\begin{theorem}
Пусть   $f\in  ZR(0,b]$. Для любых операторов $A, B\in V_{(0,b]}(X)$  таких,  что $D(A)\supseteq D(B)$ и $A-B$ принадлежит $\mathcal{I}$, оператор $f(A)-f(B)$ тоже принадлежит $\mathcal{I}$, и выполняется неравенство
$$
\|f(A)-f(B)\|_\mathcal{I}\leq M_AM_Bf'(-0)\|A-B\|_\mathcal{I}.
$$
\end{theorem}

 %\begin{proof} 
Доказательство. Пусть $A, B \in V_{(0,b]}(X)$ таковы что $D(A)\supseteq D(B)$ и $A-B\in\mathcal{I}$.
 Как показано в доказательстве теоремы 1,
\[
AR(t,A)-BR(t,B)= t R(t,A)(A-B)R(t,B), t\in (0,b], \eqno(5)
\]

Следовательно, этот оператор принадлежит  $\mathcal{I}$, и из (4) и (5) вытекает, что
\[
\| f(A)-f(B)\|_\mathcal{I}\leq M_AM_B\| A-B\|_\mathcal{I}\int\limits_0^b \frac{d\tau(t)}{t}=M_AM_Bf'(-0)\| A-B\|_\mathcal{I},
\]
поскольку  $\int_0^bd\tau(t)/t=\lim_{z\to -0}f(z)/z=f'(-0)$.
%\end{proof}

Рассматриваемое функциональное исчисление обладает следующим \textit{свойством устойчивости}.

\begin{corollary}
Если  $A_n, B_n\in V_{(0,b]}(X),\ D(A_n)\supseteq D(B_n)$, причем  $A_n-B_n\in\mathcal{I}$ и $\|A_n-B_n\|_\mathcal{I}\to 0,\ M_{A_n},  M_{B_n}\leq{\rm const}$, то $\|f(A_n)-f(B_n)\|_\mathcal{I}\to 0\ (n\to\infty)$.
\end{corollary}

\begin{remark}
Данное функциональное исчисление обладает также следующим \textit{свойством непрерывности}. Если $f_n\in ZR(0,b],\ A\in V_{(0,b]}(X)$ и  $f_n'(-0)\to 0$, то $\|f_n(A)\|\to 0$  $(n\to\infty)$.
\end{remark}

 В самом деле,  $\|f_n(A)\|\leq (1+M_A)\tau_n([0,b])\to 0$  ($\tau_n$ --- представляющая мера для $f_n$),
поскольку
$$
\frac{1}{b}\tau_n([0,b])\leq \int\limits_0^b\frac{d\tau_n(t)}{t}=\lim_{z\to -0}\frac{f_n(z)}{z}=f_n'(-0).
$$

Теорема 3 может быть применена также для оценки норм коммутаторов. Ниже через $[A,B]$ обозначается коммутатор $AB-BA$ операторов  $A$ и $B$. Если этот оператор допускает продолжение на все  $X$, то это продолжение тоже обозначается  $[A,B]$.

 \begin{corollary}
 Пусть $U$ есть автоморфизм пространства $X$, $f\in  ZR(0,b]$, $A\in V_{(0,b]}(X)$.  Если $[A,U]\in \mathcal{I}$, то $[f(A),U]\in \mathcal{I}$ и
$$
\|[f(A),U]\|_\mathcal{I}\leq M_A^2f'(-0)\|[A,U]\|_\mathcal{I}.
$$
\end{corollary}

 %\begin{proof} 
 Доказательство. Заметим, что $R(t,UAU^{-1})=UR(t,A)U^{-1}\ (t\in (0,b])$. Отсюда следует, прежде всего, что $UAU^{-1}\in V_{(0,b]}(X)$ и что $M_{UAU^{-1}}=M_A$. Более того, $f(UAU^{-1})=Uf(A)U^{-1}$, а потому $[f(A),U]=(f(A)-f(UAU^{-1}))U$. При этом   $A-UAU^{-1}=[A,U]U^{-1}\in \mathcal{I}$. Тогда по теореме 3
$$
\|[f(A),U]\|_\mathcal{I}\leq \|f(A)-f(UAU^{-1})\|_\mathcal{I}\leq M_A^2f'(-0)\|A-UAU^{-1}\|_\mathcal{I}\leq M_A^2f'(-0)\|[A,U]\|_\mathcal{I}.
$$
%\end{proof}

\begin{remark}
 Как уже отмечалось во введении,  коммутаторная липшицевость в контексте банаховых пространств являлась основным предметом пионерской работы \cite{RST}, где использовалось функциональное исчисление, основанное на развитой там технике двойных операторных интегралов для операторов скалярного типа. Поскольку операторы класса $V_{(0,b]}(X)$ не являются, вообще говоря,  операторами скалярного типа, теорема 3 и следствие 4 существенно отличаются от результатов работы \cite{RST} как по рассматриваемым классам операторов, так и по используемым классам функций.
\end{remark}

\section{Операторная дифференцируемость}

В соответствии с  \cite{mirotin_AM}, обозначим через $V_{[a,b]}(X)$ множество всех замкнутых плотно определенных операторов в $X$, для которых $[a,b]\subset\rho(A)$. Если $A\in V_{[a,b]}(X),\ f(z)=zg(z)$, где функция $g\in R[a,b]$ имеет представляющую меру  $\tau$, то мы полагаем
$$
f(A):=\int_a^bAR(t,A)d\tau(t).
$$
Функция $f(z)=zg(z)$ голоморфна в окрестности спектра оператора $A\in V_{[a,b]}(X)$ и в бесконечно удаленной точке, а потому определено ее значение на операторе $A$ в смысле голоморфного функционального исчисления. С помощью теоремы Фубини легко показать, что  функциональное исчисление, определенное выше,  является частью голоморфного.

В связи с нижеследующим  отметим, что если $A\in V_{[a,b]}(X)$, то и $A+\Delta A\in V_{[a,b]}(X)$, если $\Delta A\in \mathcal{I},\ \|\Delta A\|_\mathcal{I}<\delta_A$, где
$$
\delta_A:=\min_{\zeta\in[a,b]}\|R(\zeta,A)\|^{-1}.
$$
 Это следует из \cite[замечание IV.3.2]{Kato}.

 \begin{definition}
 Пусть $f\in ZR[a,b],\ A\in V_{[a,b]}(X)$ и пусть $\mathcal{I}$ --- операторный идеал. Ограниченный оператор $f^\nabla_A$  на $\mathcal{I}$ (трансформатор) называется $\mathcal{I}$-\textit{производной Фреше операторной функции $f$ в точке} $A$, если для  $\Delta A\in \mathcal{I}$ справедливо асимптотическое равенство
$$
\|f(A+\Delta A)-f(A)-f_A^\nabla(\Delta A)\|_\mathcal{I}=o(\|\Delta A\|_\mathcal{I}) \mbox{ при }  \|\Delta A\|_\mathcal{I}\to 0. \eqno(7)
$$
\end{definition}

   \begin{theorem}
   Функция  $f\in ZR[a,b]$ является $\mathcal{I}$-дифференцируемой по Фреше, причем ее  $\mathcal{I}$-производная Фреше в точке $A\in V_{[a,b]}(X)$ имеет вид
$$
f_A^\nabla(B)=\int\limits_a^bR(t,A)BR(t,A)td\tau(t)\ (B\in \mathcal{I}). \eqno(8)
$$
\end{theorem}

  %\begin{proof} 
  Доказательство. Прежде всего, заметим, что трансформатор $F(B)$, определяемый правой частью формулы (8), действительно ограничен, поскольку $\|F(B)\|_\mathcal{I}\leq m_A^2 \int_a^btd\tau(t)\|B\|_\mathcal{I}$, где
$$
m_A:=\max_{t\in[a,b]}\|R(t,A)\|=1/\delta_A.
$$

Далее, рассуждая, как в доказательстве теоремы 1, получаем для  $\Delta A\in \mathcal{I}$, таких, что  $A+\Delta A\in V_{[a,b]}(X)$,   равенство
$$
f(A+\Delta A)-f(A)-F(\Delta A)=\int\limits_a^b(R(t,A+\Delta A)-R(t,A))\Delta AR(t,A)td\tau(t).\eqno(9)
$$
Известно (см., например, \cite[теорема IV.1.16]{Kato}), что при $\|\Delta A\|\|R(t,A)\|<1$
$$
\|R(t,A+\Delta A)-R(t,A)\|\leq\frac{\|\Delta A\|\|R(t,A)\|^2}{1-\|\Delta A\|\|R(t,A)\|}.\eqno(10)
$$

Выберем $\Delta A$ таким, что $\|\Delta A\|_\mathcal{I}\leq 1/2m_A$. Тогда в силу (10)
$$
\|R(t,A+\Delta A)-R(t,A)\|\leq\frac{\|\Delta A\|_\mathcal{I}\|R(t,A)\|^2}{1-\|\Delta A\|_\mathcal{I}\|R(t,A)\|}\leq 2m_A^2\|\Delta A\|_\mathcal{I}.
$$
Отсюда и из равенства (9) следует, что
$$
\|f(A+\Delta A)-f(A)-F(\Delta A)\|_\mathcal{I}=o(\|\Delta A\|_\mathcal{I}) \mbox{ при }  \|\Delta A\|_\mathcal{I}\to 0,
$$
что завершает доказательство.
%\end{proof}

Для формулировки  следствия теоремы 4 введем на $V_{[a,b]}(X)$ отношение эквивалентности, считая операторы $A$ и $A'$ из $V_{[a,b]}(X)$ $\mathcal{I}$-\textit{эквивалентными}, если $A'-A\in \mathcal{I}$. Очевидно, что формула ${\rm dist}(A,A')=\|A'-A\|_\mathcal{I}$ задает метрику в каждом классе эквивалентности.

 \begin{corollary} Отображение $A\mapsto f^\nabla_A$  непрерывно в каждом классе $\mathcal{I}$-эквивалентности  операторов из $V_{[a,b]}(X)$.
 \end{corollary}

 %\begin{proof} 
 Доказательство. Пусть операторы $A$ и $A'$ из $V_{[a,b]}(X)$ эквивалентны. По теореме 4 для любого $B\in \mathcal{I}$ имеем
$$
(f^\nabla_{A'}-f^\nabla_A)B=\int\limits_a^bR(t,A')B(R(t,A')-R(t,A))td\tau(t)+\int\limits_a^b(R(t,A')-R(t,A))BR(t,A)td\tau(t).
$$
Поэтому
$$
\|(f^\nabla_{A'}-f^\nabla_A)B\|_\mathcal{I}\leq 2m_{A'}m_A\|B\|_\mathcal{I}\int\limits_a^b\|R(t,A')-R(t,A)\|td\tau(t).
$$
Пусть $\varepsilon>0$. Если оператор $\Delta A:=A'-A$ таков, что $\|\Delta A\|_\mathcal{I}<\min\{1/2m_A, \varepsilon\}$, то из формулы (10) следует, что $\|R(t,A')-R(t,A)\|\leq m_A$, а потому $m_{A'}\leq 2m_A$. Из формулы (10) следует также, что $\|R(t,A')-R(t,A)\|\leq 2m_A^2\varepsilon$. Значит, если ${\rm dist}(A,A')<\min\{1/2m_A, \varepsilon\}$, то
$$
\|(f^\nabla_{A'}-f^\nabla_A)B\|_\mathcal{I}\leq \left(8m_A^4\int\limits_a^btd\tau(t)\right)\varepsilon \|B\|_\mathcal{I},
$$
т. е.
$$
\|f^\nabla_{A'}-f^\nabla_A\|_{LB(\mathcal{I})}\leq \left(8m_A^4\int\limits_a^btd\tau(t)\right)\varepsilon,
$$
что и требовалось доказать.
%\end{proof}

Для формулировки еще одного следствия теоремы 4  заметим, что для производной функции $f\in ZR[a,b]$ с представляющей мерой $\tau$ справедливо равенство
$$
f'(z)=\int\limits_a^b\frac{td\tau(t)}{(t-z)^2}.
$$
Для функции  $f\in ZR[a,b]$ и оператора $A\in V_{[a,b]}(X)$ положим
$$
f'(A):=\int\limits_a^bR(t,A)^2td\tau(t).
$$
Ясно, что оператор $f'(A)$ ограничен.

\begin{corollary}
 Если  операторы $A$ и $B$  коммутируют, то  $f_A^\nabla(B)=f'(A)B$.
\end{corollary}

  \begin{theorem}
  Если  $f\in ZR[a,b]$, $A\in V_{[a,b]}(X)$, $B\in \mathcal{I}$, то $\mathcal{I}$-значная функция  $z\mapsto f(A+zB)-f(A)$ аналитична в  окрестности нуля $\mathcal{O}_{A, B}:=\{z\in \mathbb{\mathbb{C}}:|z|<\delta_A/\|B\|\}$, и в этой окрестности имеет место разложение
$$
f(A+zB)-f(A)=\sum\limits_{n=1}^\infty z^n C_n,\eqno(11)
$$
 где
 $$
 C_n=\frac{1}{n!}\left.\frac{d^n}{dz^n}f(A+zB)\right|_{z=0}=\int\limits_a^b(R(t,A)B)^nR(t,A)td\tau(t)\eqno(12)
 $$
(производные понимаются в смысле нормы $\|\cdot\|_\mathcal{I}$).
\end{theorem}

 %\begin{proof}
  Доказательство. Если  $z\in \mathcal{O}_{A, B}$, то, как  было отмечено перед определением 5,  $A+zB\in V_{[a,b]}(X)$. Рассуждая, как при доказательстве теоремы 3, получим, что $f(A+zB)-f(A)\in \mathcal{I}$. Пусть  $h$  настолько мало, что $z+h\in \mathcal{O}_{A, B}$. Тогда,
заменяя в формуле (7) $A$ на  $A+zB$, а $\Delta A$  --- на  $h B$ и разделив на  $h$, получим, что существует
$$
\frac{d}{dz}f(A+zB)=f^\nabla_{A+zB}(B)=\int\limits_a^bR(t,A+zB)BR(t,A+zB)td\tau(t)\eqno(13)
$$
(мы воспользовались теоремой 4). Следовательно, функция  $f(A+zB)-f(A)$ аналитична в   $\mathcal{O}_{A, B}$,
и имеет место разложение (11), в котором $C_n$ определяется первым из равенств в (12).
 Второе равенство в (12) будем доказывать по индукции. При $n=1$  оно верно в силу (13). Если предположить, что оно верно для некоторого $n$, то
$$
\left.\frac{d^{n+1}}{dz^{n+1}}f(A+zB)\right|_{z=0}=\lim\limits_{h\to 0}\frac{1}{h}\left(\left.\frac{d^{n}}{dz^{n}}f(A+zB)\right|_{z=h}-\left.\frac{d^{n}}{dz^{n}}f(A+zB)\right|_{z=0}\right)=
$$
$$
n!\lim\limits_{h\to 0}\frac{1}{h}\int\limits_a^b((R(t,A+hB)B)^nR(t,A+hB)-(R(t,A)B)^nR(t,A))td\tau(t)=
$$
$$
n!\lim\limits_{h\to 0}\frac{1}{h}\int\limits_a^b(R(t,A+hB)B)^n(R(t,A+hB)-R(t,A))td\tau(t)+
$$
$$
n!\lim\limits_{h\to 0}\frac{1}{h}\int\limits_a^b((R(t,A+hB)B)^n-(R(t,A)B)^n)R(t,A))td\tau(t)=L_1+L_2.\eqno(14)
$$
Заметим, что
$$
\|(R(t,A+hB)B)^n(R(t,A+hB)-R(t,A)\|_\mathcal{I}\leq 
$$
$$\|R(t,A+hB)\|^n\|B\|_\mathcal{I}^n\|R(t,A+hB)-R(t,A)\|.
$$
При этом из (10) следует, что при $|h|<1/(2\|B\|_\mathcal{I}m_A)$ выполняется неравенство   $\|R(t,A+hB)-R(t,A)\|<1$, а потому   $\|R(t,A+hB)\|<1+m_A$. Следовательно,
$$
\|(R(t,A+hB)B)^n(R(t,A+hB)-R(t,A)\|_\mathcal{I}\leq (1+m_A)^n\|B\|_\mathcal{I}^n.
$$
Из (10) следует также, что $\|(R(t,A+hB)-R(t,A)\|_\mathcal{I}\to 0$  при  $h\to 0$. Отсюда в силу второго резольвентного тождества выводим, что
$$
\lim\limits_{h\to 0}\left\|\frac{1}{h}(R(t,A+hB)-R(t,A))-R(t,A)BR(t,A)\right\|_\mathcal{I}= 0. \eqno(15)
$$
Воспользовавшись теоремой Лебега для интеграла Бохнера, получаем теперь, что
$$
L_1=n!\lim\limits_{h\to 0}\frac{1}{h}\int\limits_a^b(R(t,A+hB)B)^n(R(t,A+hB)-R(t,A))td\tau(t)=
$$
$$
n!\int\limits_a^b(R(t,A)B)^{n+1}R(t,A)td\tau(t).\eqno(16)
$$

При вычисления предела $L_2$ будем исходить из тождества
$$
S^n-T^n=\sum\limits_{k=1}^nS^{n-k}(S-T)T^{k-1},
$$
справедливого для любых ограниченных операторов  $S$ и $T$  в  $X$. Полагая в этом тождестве $S=R(t,A+hB)B, T=R(t,A)B$ и разделив обе его части на  $h$, заключаем, что
$$
\lim\limits_{h\to 0}\left\|\frac{1}{h}((R(t,A+hB)B)^n-(R(t,A)B)^n)-n(R(t,A)B)^{n+1}\right\|_\mathcal{I}= 0.
$$
Отсюда следует, что
$$
L_2=n!\lim\limits_{h\to 0}\frac{1}{h}\int\limits_a^b((R(t,A+hB)B)^n-(R(t,A)B)^n)R(t,A))td\tau(t)=
$$
$$
n!n\int\limits_a^b(R(t,A)B)^{n+1}R(t,A)td\tau(t).\eqno(17)
$$
Подставляя (16) и (17) в (14), получаем требуемое равенство. Теорема доказана.
%\end{proof}

\section{ Формула следа Лифшица-Крейна}

В связи со следующей теоремой отметим, что формула  следа  для  ядерных возмущений самосопряженных операторов была доказана в случае конечномерных возмущений в  \cite{L}, а в общем случае --- в \cite{K1}, \cite{K2}. Обзор более поздних исследований (в контексте гильбертовых пространств) и библиографию см. в \cite{BY}. Из последних работ на эту тему см.  \cite{Pel16}. Ниже мы вводим \textit{функцию спектрального сдвига} $\xi$ для пары операторов класса $V_{[a,b]}(X)$ и доказываем для рассматриваемого нами функционального исчисления аналог формулы следа Лифшица-Крейна. Напомним, что для банахова пространства $X$ со свойством аппроксимации существует непрерывный линейный функционал  ${\rm tr}$ нормы  1 (\textit{след}) на операторном идеале $(\mathbf{S}_1, \|\cdot \|_{\mathbf{S}_1})$ ядерных операторов, определенных на  $X$ (см., например, \cite[с. 64]{DF}).

 \begin{theorem}
 Пусть пространство $X$ обладает свойством аппроксимации, $A, B\in V_{[a,b]}(X)$ $D(A)\supseteq D(B)$,  $A-B\in \mathbf{S}_1$. Тогда
 существует  функция $\xi$, аналитическая в окрестности отрезка $[a,b]$ и такая, что
 для любой функции $f\in ZR[a,b]$ справедливо равенство
$$
{\rm tr}(f(A)-f(B))=\frac{1}{2\pi i}\int\limits_\Gamma \xi(z)f'(z)dz,
$$
где контур $\Gamma$ охватывает $[a,b]$ и лежит в связной компоненте множества  $\rho(A)\cap\rho(B)$, содержащей $[a,b]$.
\end{theorem}

 %\begin{proof} 
 Доказательство. Тот факт, что $f(A)-f(B)\in \mathbf{S}_1$, доказывается так же, как и в теореме 3. При этом в силу аналога формулы  (4) (в которой нижний предел интегрирования заменен на $a$) и формулы (5) справедливо равенство
$$
 f(A)-f(B)=\int\limits_a^b R(t,A)(A-B)R(t,B)td\tau(t),
 $$
причем интеграл в правой части  существует в смысле Бохнера  относительно ядерной нормы. Следовательно,
$$
{\rm tr}(f(A)-f(B))=\int\limits_a^b\varphi(t)td\tau(t),
$$
где $\varphi(t)={\rm tr}(R(t,A)(A-B)R(t,B))$. Последняя функция аналитична в односвязной окрестности отрезка $[a,b]$, содержащейся в $\rho(A)\cap\rho(B)$. В самом деле, поскольку отрезок $[a,b]$ лежит в $\rho(A)\cap\rho(B)$, для некоторой окрестности любой точки $z_0\in [a,b]$ справедливы абсолютно сходящиеся по норме оператора разложения $R(z,A)=\sum_{n=0}^\infty (z-z_0)^nA_n$, $R(z,B)=\sum_{m=0}^\infty (z-z_0)^mB_m$ с ограниченными операторными коэффициентами. Поэтому
$$
R(z,A)(A-B)R(z,B)=\sum_{n=0}^\infty\sum_{m=0}^\infty (z-z_0)^{n+m}A_n(A-B)B_m,
$$
причем ряд, стоящий в правой части, абсолютно сходится в ядерной норме,
поскольку
$$
\|(z-z_0)^{n+m}A_n(A-B)B_m\|_{\mathbf{S}_1}\leq |z-z_0|^{n+m}\|A_n\|\|B_m\|\|(A-B)\|_{\mathbf{S}_1}.
$$
Если $\xi$ есть первообразная  функции $\varphi$ в достаточно малой односвязной окрестности отрезка  $[a,b]$ (первообразная существует, так как  $\varphi$ аналитична в такой окрестности), то, как известно, по формуле Коши
$$
\varphi(t)=\frac{1}{2\pi i}\int\limits_\Gamma \frac{\xi(z)}{(z-t)^2}dz,
$$
где контур $\Gamma$ охватывает $[a,b]$ и лежит в выбранной окрестности.

Тогда, в силу теоремы Фубини,
$$
{\rm tr}(f(A)-f(B))=\frac{1}{2\pi i}\int\limits_\Gamma \xi(z)\int\limits_a^b\frac{td\tau(t)}{(z-t)^2}dz=\frac{1}{2\pi i}\int\limits_\Gamma \xi(z)f'(z)dz
$$
(теорема Фубини применима, поскольку подынтегральная функция в соответствующем двойном интеграле непрерывна на компакте $\Gamma\times[a,b]$), что и требовалось доказать.
%\end{proof}

% Литература

\end{document}